\documentclass{amsart}

\usepackage[frame,cmtip,arrow,matrix,line,graph,curve]{xy}
\usepackage{amsmath,amssymb, amscd}
\usepackage{latexsym}
\usepackage{graphics}

\numberwithin{equation}{section}

\newtheorem{theorem}{Theorem}[section]
\newtheorem{proposition}[theorem]{Proposition}
\newtheorem{corollary}[theorem]{Corollary}
\newtheorem{lemma}[theorem]{Lemma}

\newtheorem{remit}[theorem]{Remark}
\newtheorem{definit}[theorem]{Definition}
\newtheorem{question}[theorem]{Question}

\newenvironment{definition}{\begin{definit}\rm}{\end{definit}}

\newenvironment{remark}{\begin{remit}\rm}{\end{remit}}

\newcommand{\pp}{\mathbb{P}}
\newcommand{\qq}{\mathbb{Q}}
\newcommand{\cc}{\mathbb{C}}

\newcommand{\cO}{\mathcal{O} }

\begin{document}

\title[Symplectic desingularization of moduli space of sheaves]{Symplectic
desingularization of moduli space of sheaves on a K3 surface}
\date{}

\author{Young-Hoon Kiem}
\address{Department of Mathematics, Seoul National University, Seoul
151-747, Korea} \email{kiem@math.snu.ac.kr}

\thanks{Young-Hoon Kiem was partially supported by KOSEF and KRF}


 \maketitle


\section{Introduction}
Let $X$ be a projective K3 surface with generic polarization
$\cO_X(1)$ and let $M_c=M(2,0,c)$ be the moduli space of
semistable torsion-free sheaves on $X$ of rank $2$, with Chern
classes $c_1=0$ and $c_2=c$. When $c=2n\ge 4$ is even, $M_c$ is a
singular projective variety. Recently O'Grady raised the following
question (\cite{ogrady} 0.1).
\begin{question}\label{q1.1} Does there exist a symplectic desingularization
of $M_{2n}$?\end{question}

In \cite{ogrady}, he analyzes Kirwan's desingularization
$\widehat{M}_c$ of $M_c$ and proves that $\widehat{M}_c$ can be
blown down twice and that as a result he gets a symplectic
desingularization $\widetilde{M}_c$ of $M_c$ in the case when
$c=4$. This turns out to be a new irreducible symplectic variety.

When $c\ge 6$, O'Grady conjectures that there is no smooth
symplectic model of $M_c$ (\cite{ogrady} page 50). The purpose of
this paper is to provide a partial answer to Question \ref{q1.1}.
\begin{theorem}\label{th1.2} There is no symplectic desingularization of
$M_{2n}$ if $\frac{n\,a_n}{2n-3}$ is not an integer where $a_n$ is
the Euler number of the Hilbert scheme $X^{[n]}$ of $n$ points in
$X$.\end{theorem}

It is well-known that $a_n$ is given by the equation
$$\sum_{n=0}^\infty a_nq^n=\prod_{n=1}^{\infty}1/(1-q^n)^{24}.$$
By direct computation, one can check that $\frac{n\,a_n}{2n-3}$ is
not an integer for $n=$ 5, 6, 8, 11, 12, 13, 15, 16, 17, 18, 19,
20,$\cdots$.

The idea of the proof is to use motivic integral in the form of
stringy Euler numbers. If there is an irreducible symplectic
desingularization $\widetilde{M}_c$ of $M_c$, then the stringy
Euler number of $M_c$ is equal to the ordinary Euler number of
$\widetilde{M}_c$ because the canonical divisors of both
$\widetilde{M}_c$ and $M_c$ are trivial (Theorem \ref{crepthm}).
In particular, we deduce that the stringy Euler number
$e_{st}(M_c)$ must be an integer. Therefore, Theorem \ref{th1.2}
is a consequence of the following.
\begin{proposition} \label{prop1} The stringy Euler number $e_{st}(M_{2n})$ is
of the form $$\frac{n\,a_n}{2n-3} +\text{ integer }.$$
\end{proposition}
We prove this proposition in section \ref{mainproof} after a brief
review of preliminaries.

One motivation for Question \ref{q1.1} is to find a mathematical
interpretation of Vafa-Witten's formula (\cite{VW} 4.17) which
says the ``Euler characteristic'' of $M_{2n}$ is
$$e^{VW}(M_{2n})=a_{4n-3}+\frac14a_n.$$ Because $k/4\neq l/(2n-3)$
for $1\le k\le 3$, $1\le l< 2n-3$, we deduce the following from
Proposition \ref{prop1}.
\begin{corollary} The stringy Euler number $e_{st}(M_{2n})$
is not  Vafa-Witten's Euler characteristic $e^{VW}(M_{2n})$ in
general.\end{corollary}

It is my great pleasure to express
my gratitude to Professor Jun Li for useful discussions and to
Jaeyoo Choy for motivating questions.

\section{Preliminaries}
In this section, we recall the definition and basic facts about
stringy Euler numbers. The references are \cite{Bat, DenLoe1}.

Let $W$ be a variety with at worst \emph{log-terminal}
singularities, i.e. \begin{itemize}
\item  $W$ is $\qq$-Gorenstein
\item for a resolution of singularities $\rho:V\to W$ such that
the exceptional locus of $\rho$ is a divisor $D$ whose irreducible
 components $D_1, \cdots, D_r$ are smooth divisors with only normal
  crossings, we have
$$K_V=\rho^*K_W+\sum_{i=1}^ra_iD_i$$
with $a_i>-1$ for all $i$, where $D_i$ runs over all irreducible
components of $D$. The divisor $\sum_{i=1}^ra_iD_i$ is called the
\emph{discrepancy divisor}.\end{itemize}

For each subset $J\subset I=\{1,2,\cdots,r\}$, define
$D_J=\cap_{j\in J}D_j$, $D_{\emptyset}=Y$ and
$D_J^{0}=D_J-\cup_{j\in I-J}D_j$. Then the stringy E-function of
$W$ is defined by
\begin{equation}\label{EstDef}
E_{st}(W;u,v)=\sum_{J\subset I}E(D_J^0;u,v)\prod_{j\in
J}\frac{uv-1}{(uv)^{a_j+1}-1}\end{equation} where
$$E(Z;u,v)=\sum_{p,q}\sum_{k\ge 0} (-1)^kh^{p,q}(H^k_c(Z;\cc))u^pv^q$$
 is the Hodge-Deligne polynomial for a variety $Z$. Note that the
 Hodge-Deligne polynomials have\begin{itemize}
\item  the additive property:
 $E(Z;u,v)=E(U;u,v)+E(Z-U;u,v)$ if $U$ is a smooth open subvariety of $Z$
 \item the multiplicative property:
$E(Z;u,v)=E(B;u,v)\, E(F;u,v)$ if $Z$ is a locally trivial
$F$-bundle over $B$.\end{itemize}

\begin{definition} The stringy Euler number is defined as
\begin{equation}\label{est} e_{st}(W)=\lim_{u,v\to 1}E_{st}(W;u,v)= \sum_{J\subset
I}e(D_J^0)\prod_{j\in J}\frac{1}{a_j+1}\end{equation} where
$e(D_J^0)=E(D_J^0;1,1)$.
\end{definition}

The ``change of variable formula" (Theorem 6.27 in \cite{Bat},
Lemma 3.3 in \cite{DenLoe1}) implies that the function $E_{st}$ is
independent of the choice of a resolution and the following holds.

\begin{theorem}\label{crepthm} {\rm (\cite{Bat} Theorem 3.12)}
Suppose $W$ is a $\qq$-Gorenstein algebraic variety with at worst
log-terminal singularities. If $\rho:V\to W$ is a crepant
desingularization (i.e. $\rho^*K_W=K_V$) then
$E_{st}(W;u,v)=E(V;u,v)$. In particular, $e_{st}(W)=e(V)$ is an
integer.
\end{theorem}

\section{Proof of Proposition \ref{prop1}}\label{mainproof}

We fix a generic polarization of $X$ as in \cite{ogrady} page 50.
The moduli space $M_{2n}$ has a stratification
$$M_{2n}=M^s_{2n}\sqcup (\Sigma-\Omega) \sqcup \Omega$$
where $M^s_{2n}$ is the locus of stable sheaves and $$\Sigma\cong
(X^{[n]}\times X^{[n]})/\text{involution} $$ is the locus of
sheaves of the form $I_{Z}\oplus I_{Z'}$ ($[Z], [Z']\in X^{[n]}$)
while $$\Omega\cong X^{[n]}$$ is the locus of sheaves $I_Z\oplus
I_Z$. Kirwan's desingularization $\rho:\widehat{M}_{2n}\to M_{2n}$
is obtained by blowing up $M_c$ first along $\Omega$, next along
the proper transform of $\Sigma$ and finally along the proper
transform of a subvariety $\Delta$ in the exceptional divisor of
the first blow-up. This is indeed a desingularization by
\cite{ogrady} Proposition 1.8.3.

Let $D_1=\widehat{\Omega}$, $D_2=\widehat{\Sigma}$,
$D_3=\widehat{\Delta}$ be the (proper transforms of the)
exceptional divisors of the three blow-ups. Then they are smooth
divisors with only normal crossings and the discrepancy divisor of
$\rho:\widehat{M}_{2n}\to M_{2n}$ is (\cite{ogrady} 6.1)
$$(6n-7)D_1+(2n-4)D_2+(4n-6)D_3 .$$
Therefore the singularities are terminal for $n\ge 2$ and from
\eqref{est} the stringy Euler number of $M_{2n}$ is given by
\begin{equation}\label{eq3.1}\begin{array}{ll}
&e(M^s_{2n})+e(D_1^0)\frac1{6n-6}+e(D_2^0)\frac1{2n-3}+e(D_3^0)\frac1{4n-5}\\
&+e(D_{12}^0)\frac1{6n-6}\frac1{2n-3}
+e(D_{23}^0)\frac1{2n-3}\frac1{4n-5}\\
&+e(D_{13}^0)\frac1{6n-6}\frac1{4n-5}
+e(D_{123}^0)\frac1{6n-6}\frac1{2n-3}\frac1{4n-5}\quad
.\end{array}
\end{equation}

We need to compute the (virtual) Euler numbers of $D_J^0$ for
$J\subset \{1,2,3\}$. Let $(E,\omega)$ be a symplectic vector
space of dimension $c=2n$. Let $\mathrm{Gr}^{\omega}(k,c)$ be the
Grassmannian of $k$ dimensional subspaces of $E$ isotropic with
respect to the symplectic form $\omega$ (i.e. the restriction of
$\omega$ to the subspace is zero).
\begin{lemma}\label{lem3.1} For $k\le n$, the Euler number of $\mathrm{Gr}^{\omega}(k,2n)$ is
$2^k\binom{n}{k}$.\end{lemma}
\begin{proof} Consider the incidence variety $$\{(a,b)\in
\mathrm{Gr}^{\omega}(k-1,2n)\times \mathrm{Gr}^{\omega}(k,2n)\,|\,
a\subset b\}.$$ This is a $\pp^{2n-2k+1}$-bundle over
$\mathrm{Gr}^{\omega}(k-1,2n)$ and a $\pp^{k-1}$-bundle over
$\mathrm{Gr}^{\omega}(k,2n)$. The formula follows from an
induction on $k$.\end{proof}

Let $\hat\pp^5$ be the blow-up of $\pp^5$ (projectivization of the
space of $3\times 3$ symmetric matrices) along $\pp^2$ (the locus
of rank 1 matrices). We have the following from \cite{ogrady} \S6
and \cite{ogrady1} \S3.
\begin{proposition}\label{prop3.2}\begin{enumerate}
\item $D_1$ is a $\hat\pp^5$-bundle over a
$\mathrm{Gr}^{\omega}(3,2n)$-bundle over $X^{[n]}$.
\item $D_2^0$ is a $\pp^{2n-4}$-bundle over a $\pp^{2n-3}$-bundle over
{\rm $(X^{[n]}\times X^{[n]}-X^{[n]})/\text{involution}$}.
\item $D_3$ is a $\pp^{2n-4}\times\pp^2$-bundle over
a $\mathrm{Gr}^{\omega}(2,2n)$-bundle over $X^{[n]}$.
\item $D_1\cap D_2$ is a $\pp^2\times \pp^2$-bundle over
 $\mathrm{Gr}^{\omega}(3,2n)$-bundle over $X^{[n]}$.
\item $D_2\cap D_3$ is a $\pp^{2n-4}\times\pp^1$-bundle
over a $\mathrm{Gr}^{\omega}(2,2n)$-bundle over $X^{[n]}$.
\item $D_1\cap D_3$ is a $\pp^2\times\pp^{2n-5}$-bundle
over a $\mathrm{Gr}^{\omega}(2,2n)$-bundle over $X^{[n]}$.
\item $D_1\cap D_2\cap D_3$ is a $\pp^1\times\pp^{2n-5}$-bundle
over a $\mathrm{Gr}^{\omega}(2,2n)$-bundle over $X^{[n]}$.
\end{enumerate}\end{proposition}
For instance, (1) is just Proposition 6.2 of \cite{ogrady} and (2)
is Proposition 3.3.2 of \cite{ogrady1} while (3) is Lemma 3.5.4 in
\cite{ogrady1}.

From Proposition \ref{prop3.2} and Lemma \ref{lem3.1}, we have the
following by the additive and multiplicative properties of the
(virtual) Euler numbers:
$$\begin{array}{lll}
 e(D_1^0)=0, & e(D_2^0)= (2n-3)(2n-2) \frac12 (a_n^2-a_n),\\
 e(D_3^0)=2^2\binom{n}{2}\,a_n, & e(D_{12}^0)=3\cdot
 2^3\binom{n}{3}\,a_n\\
 e(D_{23}^0)=2\cdot 2^2\binom{n}{2}\,a_n, &
 e(D_{13}^0)=(2n-4)2^2\binom{n}{2}\,a_n\\
 e(D_{123}^0)=2(2n-4)2^2\binom{n}{2}\,a_n\quad .
\end{array}$$
Hence from the formula \eqref{eq3.1}, the stringy Euler number of
$M_{2n}$ is given by
$$e_{st}(M_{2n})=e(M^s_{2n})+(n-1)(a_n^2-a_n)+n\frac{2n-2}{2n-3}a_n
=\frac{n\,a_n}{2n-3} + \text{ integer } $$ since $e(M^s_{2n})$ is
an integer. So we proved Proposition \ref{prop1}.

\begin{remark} For the moduli space of rank 2 bundles over a
smooth projective curve, the stringy E-function and the stringy
Euler number are computed in \cite{kiem} and \cite{KL}.
\end{remark}

\end{document}